# A New Numerical Method for Fast Solution of Partial Integro-Differential Equations

P. Dourbal[1] and M. Pekker[2]

[1] Dourbal Electric, Inc., 10 Schalk's Crossing Rd., St. 501-295, Princeton Jct., NJ 08536

*paul@dourbalelectric.com*

[2] The George Washington University, Science & Engineering Hall 3550, 800 22nd Street, Northwest Washington, DC 20052, *pekkerm@gmail.com*

## Abstract

We propose a new method of numerical solution for partial differential equations. The method is based on a fast matrix multiplication algorithm developed by Dourbal [1,2]. We use a two-dimensional Poison equation for comparison of the proposed method with conventional numerical methods. We have shown that for $N \geq 50 \times 50$ grids, the new method allows for linear growth in the number of elementary addition and multiplication operations with the growth of N, as contrasted with quadratic growth necessitated by the standard numerical methods [4 - 8]. The algorithm described here can be easily generalized for any differential equations, and is not specific to the Poison equation.

## Introduction

Many problems in physics and engineering can be described in terms of multivariate (2- or 3-dimensional or more) partial differential equations. As an example, let's consider a modified Poisson equation describing an electrical field distribution in a charged media:

$$\frac{\partial}{\partial x_i}\varepsilon_i\frac{\partial u}{\partial x_i}=\rho \quad (1)$$

where $\varepsilon_i(x_i)$ is a permittivity of the media, $\rho(x_i,t)$ is the charge density, $x_i$ is an *i*-spatial coordinate, and *t* is time. A matrix corresponding to a differential operator $\frac{\partial}{\partial x_i}\varepsilon_i\frac{\partial}{\partial x_i}$ from (1) has size $N^2$ with $N=n_x\cdot n_y\cdot n_z$, where $n_x,n_y,n_z$ are the numbers of grid points along the $x,y,z$ axes. Eq. (1) is usually solved numerically by using the iterative methods described in [3].

Let a square matrix $L_\varepsilon^{-1}$ corresponds to the inverse of the $\frac{\partial}{\partial x_i}\varepsilon_i\frac{\partial}{\partial x_i}$ operator.

At first glance, it looks preferable to multiply the matrix $L_\varepsilon^{-1}$ by vector $\vec{\rho}=\rho(x_i)$ at every time moment. However, when the number of elementary addition and multiplication operations to obtain the products of matrix $L_\varepsilon^{-1}$ by vector $\vec{\rho}$ is prohibitively high, as it proportional to $N^2$, iterative methods are used.

The proposed numerical method for solving Eq. (1) reduces the number of operations required to obtain each $L_\varepsilon^{-1}\cdot\vec{\rho}$ product by a factor of $N$. In this method, the total number of elementary multiplications and additions in each matrix multiplication for solving Eq. (1) are correspondingly $N_\times=\alpha_\times\cdot N$ and $N_+=\alpha_+\cdot N$, where coefficients $\alpha_\times$ and $\alpha_+$ grow proportionally to the precision of the inverted matrix elements. Therefore, our proposed direct method for solving Eq. (1) becomes preferable as the size of matrix $L_\varepsilon$ increases.

**Construction of a system of linear equations numerically solving Eq. (1), matrix $L_\varepsilon$**

As the first example, let us consider a standard second-order differential schema for Eq. (1) on a two dimensional grid with $\varepsilon = 1$:

$$
\begin{aligned}
& h = L/(n-1) \\
& x_i = h \cdot (i-1) \quad \mathrm{i} = 1,2,..,\mathrm{n} \\
& y_j = h \cdot (j-1) \quad \mathrm{j} = 1,2,..,\mathrm{n} \\
& u_{i,j} - 0.25 \cdot \left(u_{i-1,j} + u_{i+,j} + u_{i,j-1} + u_{i,j+1}\right) = -\rho_{i,j} \cdot h^2
\end{aligned}
\tag{2}
$$

A graphical representation of Eq. (2). is shown on Fig. 1 where dots and crosses represent correspondingly the internal and the external vertices of the grid.

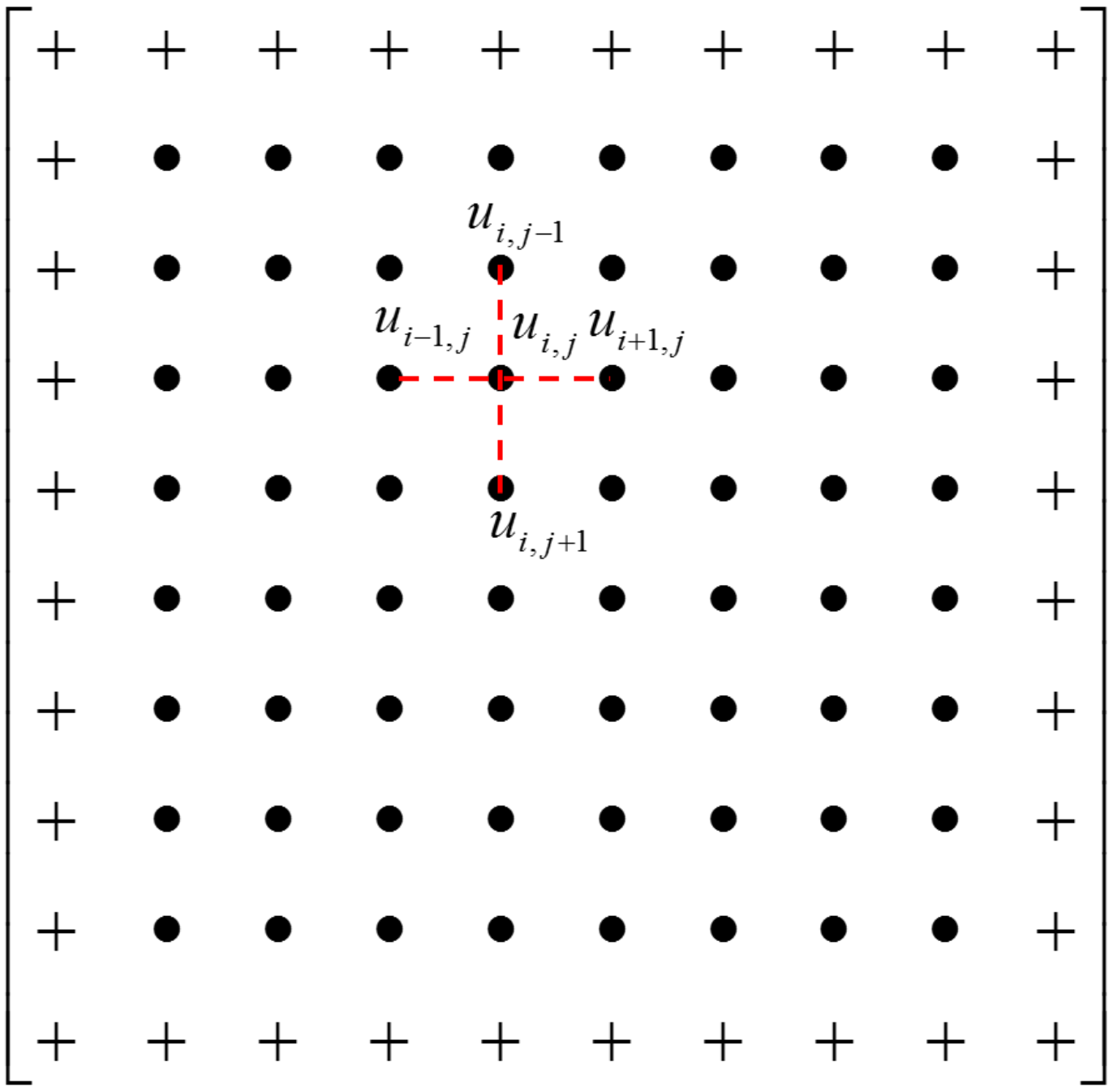

*Fig. 1. A graphical representation of the numerical differentiation operator $L_\varepsilon$. The internal vertices of the grid are represented by the dots and the external - by the crosses.*

In general, $\varepsilon$ is a function of coordinates, the grid is non-uniform, and the coefficients in front of the $u_{i,j,}$ $u_{i-1,}$ $u_{i+1,j,}$ $u_{i,j-1,}$ $u_{i,j+1}$ in Eq. (2) are not equal to unity. Those conditions do not significantly change the structure of the inverse of matrix (2) and the results presented below are true for cases with $\varepsilon = \varepsilon(x_i)$ and for non-uniform grids.

A grid for $n_x = 3, n_y = 5$ is shown on Fig. 2, and the corresponding system of equations is shown on Fig. 3. with the unknowns $u_{2,2}$ , $u_{3,2}, u_{4,2}$ and the boundary conditions

$u_{1,1}, u_{1,2}, u_{1,3}, u_{2,1}, u_{2,3}, u_{3,1}, u_{3,3}, u_{4,1}, u_{4,3}, u_{5,1}, u_{5,2}, u_{5,3}$.

$$\begin{bmatrix} + & + & + \\ + & \bullet & + \\ + & \bullet & + \\ + & \bullet & + \\ + & + & + \end{bmatrix}$$

*Fig. 2. A grid corresponding to Eq. (2) with* $n_x = 3$, $n_y = 5$. *Crosses correspond to the boundary conditions.*

$$\begin{bmatrix}
1 & 0 & 0 & 0 & 0 & 0 & 0 & 0 & 0 & 0 & 0 & 0 & 0 & 0 & 0 \\
0 & 1 & 0 & 0 & 0 & 0 & 0 & 0 & 0 & 0 & 0 & 0 & 0 & 0 & 0 \\
0 & 0 & 1 & 0 & 0 & 0 & 0 & 0 & 0 & 0 & 0 & 0 & 0 & 0 & 0 \\
0 & 0 & 0 & 1 & 0 & 0 & 0 & 0 & 0 & 0 & 0 & 0 & 0 & 0 & 0 \\
0 & -0.25 & 0 & -0.25 & 1 & -0.25 & 0 & -0.25 & 0 & 0 & 0 & 0 & 0 & 0 & 0 \\
0 & 0 & 0 & 0 & 0 & 1 & 0 & 0 & 0 & 0 & 0 & 0 & 0 & 0 & 0 \\
0 & 0 & 0 & 0 & 0 & 0 & 1 & 0 & 0 & 0 & 0 & 0 & 0 & 0 & 0 \\
0 & 0 & 0 & 0 & -0.25 & 0 & -0.25 & 1 & -0.25 & 0 & -0.25 & 0 & 0 & 0 & 0 \\
0 & 0 & 0 & 0 & 0 & 0 & 0 & 0 & 1 & 0 & 0 & 0 & 0 & 0 & 0 \\
0 & 0 & 0 & 0 & 0 & 0 & 0 & 0 & 0 & 1 & 0 & 0 & 0 & 0 & 0 \\
0 & 0 & 0 & 0 & 0 & 0 & 0 & -0.25 & 0 & -0.25 & 1 & -0.25 & 0 & -0.25 & 0 \\
0 & 0 & 0 & 0 & 0 & 0 & 0 & 0 & 0 & 0 & 0 & 1 & 0 & 0 & 0 \\
0 & 0 & 0 & 0 & 0 & 0 & 0 & 0 & 0 & 0 & 0 & 0 & 1 & 0 & 0 \\
0 & 0 & 0 & 0 & 0 & 0 & 0 & 0 & 0 & 0 & 0 & 0 & 0 & 1 & 0 \\
0 & 0 & 0 & 0 & 0 & 0 & 0 & 0 & 0 & 0 & 0 & 0 & 0 & 0 & 1
\end{bmatrix} \bullet \begin{bmatrix} u_{1,1} \\ u_{1,2} \\ u_{1,3} \\ u_{2,1} \\ u_{2,2} \\ u_{2,3} \\ u_{3,1} \\ u_{3,2} \\ u_{3,3} \\ u_{4,1} \\ u_{4,2} \\ u_{4,3} \\ u_{5,1} \\ u_{5,2} \\ u_{5,3} \end{bmatrix} = \begin{bmatrix} u_{1,1} \\ u_{1,2} \\ u_{1,3} \\ u_{2,1} \\ \rho_{2,2} \\ u_{2,3} \\ u_{3,1} \\ \rho_{3,2} \\ u_{3,3} \\ u_{4,1} \\ \rho_{4,2} \\ u_{4,3} \\ u_{5,1} \\ u_{5,2} \\ u_{5,3} \end{bmatrix}$$

*Fig. 3. A system of equations corresponding to Eq. (2) for the case of* $n_x = 3, n_y = 5$. *The diagonal elements of the matrix correspond to the points where the boundary conditions are defined.*

## The solution steps

The solution process includes the following steps:

1. For a given grid size $N = n_x \times n_y \times n_z$ we generate a square $N \times N$ element matrix $[A]_{N,N}$ corresponding to the differential operator $\frac{\partial}{\partial x_i} \varepsilon_i \frac{\partial}{\partial x_i}$.
2. Using any standard matrix inversion method, we find the matrix $[A]_{N,N}^{-1}$ - the inverse of the matrix $[A]_{N,N}$.
3. We find the vector $u_N$ solving Eq. (2).

The above steps 1 and 2 use standard methods.

Fig. 4 shows how the size of the matrix relates to the total number of elementary multiplication $N_\times$ and addition $N_+$ operations required to multiply an inverted matrix $[A]_{N,N}^{-1}$ by a vector using the proposed method. The curves on Fig. 4 correspond to the number of significant figures after the decimal point in the representation of the inverse matrix $[A]_{N,N}^{-1}$. Dashed lines show an asymptotic behavior for $N_\times = \alpha_\times \cdot N$, $N_+ = \alpha_+ \cdot N$.

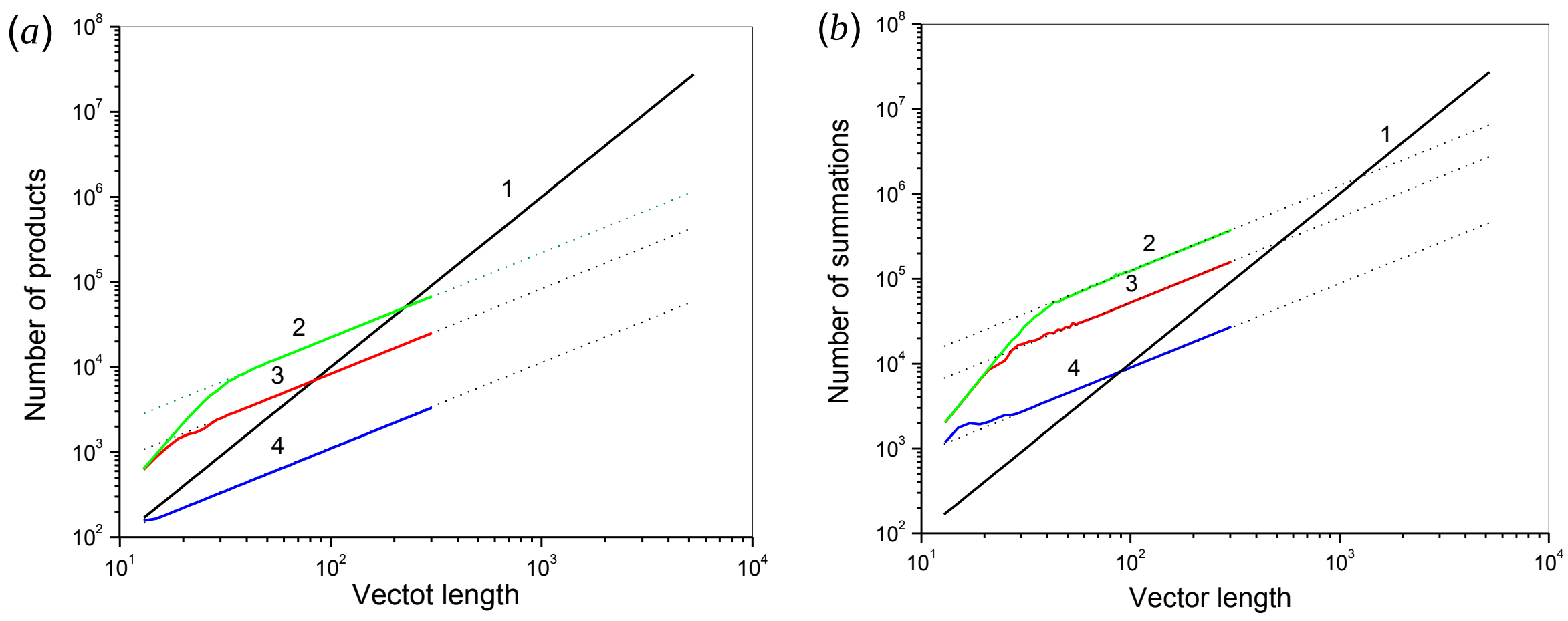

*Fig. 4. The number of elementary operations required to multiply matrix* $[A]^{-1}_{N,N}$ *by a vector as a function of the matrix size* $N \times N$ *: (a) multiplications, (b) additions. Line 1 corresponds to the conventional method of matrix by vector multiplication; line 2 corresponds to the proposed method with 6 decimal digits precision, line 3 - with 4 digits precision, and line 4 - with 2 digits precision. Dashed lines show the asymptotic behavior.*

As one can see from Fig. 4, for a given number of significant figures $m$ in the inverse matrix components, the number of elementary addition and multiplication operations required to compute the product of the matrix by a vector grows linearly with the size of the matrix. Therefore, the computational complexity of the new method is significantly reduced for $N^2 >> 10^m$.

Fig. 5 shows the relationship between the coefficients $\alpha_\times$ and $\alpha_+$ and the number of significant figures $m$ in inverse matrix $[A]^{-1}_{N,N}$.

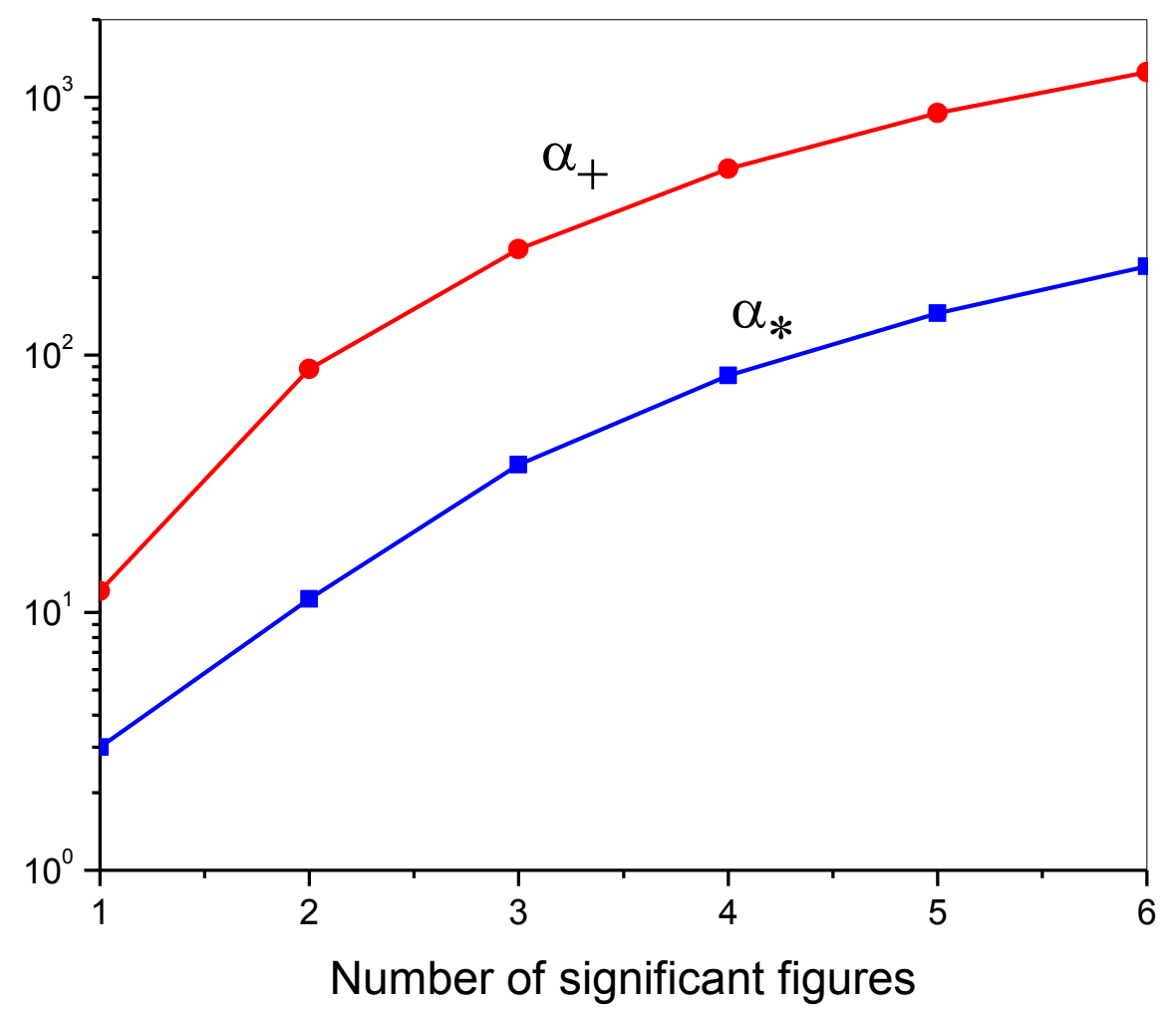

*Fig. 5. Dependence of the coefficients* $\alpha_\times$ *and* $\alpha_+$ *on the number of significant figures* $m$ *in inverse matrix* $[A]^{-1}_{N,N}$.

To show that the above results do not depend neither on $\varepsilon$ being a variable or a constant, nor on a grid having a constant or variable step size, we examined the following linear systems of differential equations:

$$u^k_{i,j} - 0.25\cdot\left(\beta_{i-1,j}u^k_{i-1,j} + \beta_{i+1,j}u^k_{i+1,j} + \beta_{i,j-1}u^k_{i,j-1} + \beta_{i,j+1}u^k_{i,j+1}\right) = \rho_{i,j} \qquad (3)$$

We have found that even for the worst case, where coefficients $\beta_{i,j}$ are arbitrary numbers, reduction in the computational intensity achieved by the use of the Dourbal method remains the same as in the above example.

## Numerical example

To demonstrate that the results do not depend on whether $\varepsilon$ is a variable or a constant, or on whether the grid is uniform or non-uniform, let us examine this equation

$$\frac{\partial^2 u}{\partial^2 x}+\frac{\partial^2 u}{\partial^2 y}=\left(12x^2-6x\right)\cdot\left(y^3-y^2\right)+\left(x^4-x^3\right)\cdot\left(6y-2\right) \tag{4}$$

defined on the domain $0 \le x \le 1,\ 0 \le y \le 1$ with zero boundary conditions. The exact solution for this equation is

$$u=\left(x^4-x^3\right)\cdot\left(y^3-y^2\right) \tag{5}$$

Following the procedure described above, we construct matrices corresponding to Eq. (4) with the number of grid points along each of the coordinate axes equal to $n_x = n_y$ = 5, 11, 21, 41, 81. The resulting matrices are inverted and then the elements of each are rounded, leaving 1, 2, 3 and 6 digits after the decimal point. Then, the vectors corresponding to the right side of Eq. (5) and the boundary conditions (an example of such a vector is given in Fig. 3) are multiplied by the resulting matrices.

Fig. 6 shows the exact analytical solution for Eq. (4).

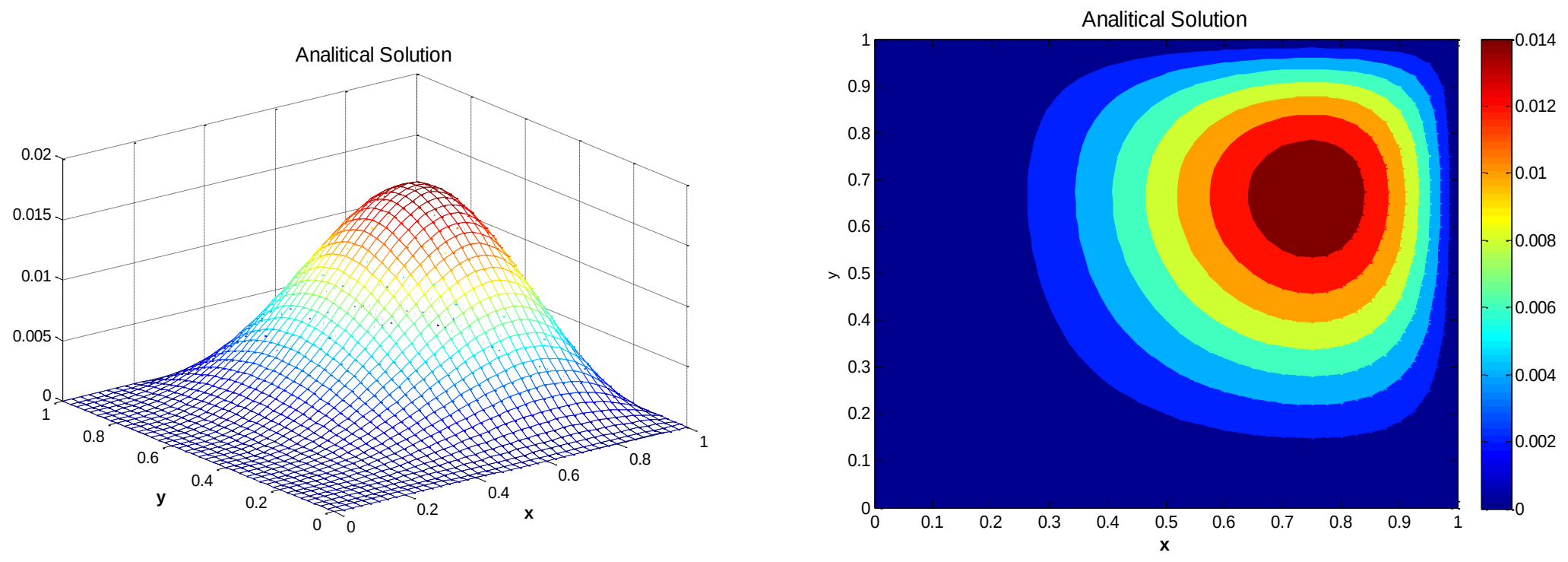

*Figure 6. Analytical solution for Eq. (4).*

Table 1 shows the error in the numerical solutions for grid sizes $n_x = n_y$ = 5, 11, 21, 41, 81 without rounding and with rounding of the inverted matrix to 1, 2, 3, and 6 decimal places. The error is computed as the maximum difference between the numerical and the exact solutions over all points of the grid, divided by the maximum value of the exact solution.

$$\varepsilon = \frac{\max\left|u_{i,j}^{anal} - u_{i,j}^{num}\right|}{\max\left|u_{i,j}^{anal}\right|} \tag{6}$$

| **Rounding** | **Error** | | | | |
|---|---|---|---|---|---|
| Digits after decimal point | $n_x$=$n_y$=5 | $n_x$=$n_y$=11 | $n_x$=$n_y$=21 | $n_x$=$n_y$=41 | $n_x$=$n_y$=81 |
| 1 | 0.0569 | 0.0352 | 0.0292 | 0.02843 | 0.02829 |
| 2 | 0.0658 | 0.0112 | 0.0038 | 0.00276 | 0.00277 |
| 3 | 0.0658 | 0.0105 | 0.0026 | 0.00074 | 0.00027 |
| 5 | 0.0658 | 0.0107 | 0.0026 | 0.00065 | 0.00016 |
| Without rounding | 0.0658 | 0.0107 | 0.0026 | 0.00065 | 0.00016 |

*Table 1. Error in the numerical solutions without and with rounding of the inverted matrix.*

.

As one can see from table 1, the relative error of the numerical solution is below 0.1 % as compared to the exact solution if the rounding is two or three digits after the decimal point for the grid size more than 250.

It must be stressed that in the last few decades several algorithms have been developed which allow for reductions in the number of multiplication and addition steps required to multiply two matrices. Pan proposes an algorithm [4] in which the multiplication of two $n \times n$ matrices requires $n^{2.81}$ multiplication and $7 \cdot n^{2.81}$ addition operations. Other algorithms were proposed in [5 - 8] with number of operations reduced by not more than the order of $\sim n^{2.37}$. These algorithms are complex and rarely used.

### Computational complexity of the Dourbal method compared to conventional iteration methods

There are two most widely used conventional iterative methods for solving linear systems of equations (2) - successive over-relaxation (SOR) and Seidel [9]. As shown in [9], in iterative methods, the number of iterations is proportional to the required solution accuracy $\varepsilon$ (6) and the width of matrix (2), i.e. $N = n_x \cdot n_y$:

$$I_{N,\varepsilon} \approx \frac{2N}{\pi} \ln\left(\frac{1}{\varepsilon}\right) = \frac{2\ln(10)}{\pi} m \cdot N \tag{7}$$

Here $\varepsilon = 10^{-m}$.

Therefore, the number of elementary mathematical operations for one iteration in (2) is $S_1 = 5N$, and the total number of elementary operations is:

$$S_{I,total} = S_1 \cdot I_{N,\varepsilon} = \frac{10\ln(10)}{\pi} m \cdot N^2 \tag{8}$$

In contrast, the Dourbal method is direct. It performs a fast multiplication of inverse matrix (3). Asymptotic number of elementary product and summation operations is:

$$S_{D,prod} = 2.5m^{2.5} \cdot N \text{ for elementary products,} \tag{9}$$

and

$S_{D,sum} = 14.2m^{2.5} \cdot N$ for elementary summations

(10)

From (8), it is evident that the number of elementary operations in conventional iterative methods grows proportionally to $N^2$. This is comparable to the conventional direct method. In the Dourbal method, the number of elementary operations is proportional to $N$, as may be seen from (9) and (10) . Therefore, the Dourbal method allows for linear, rather than quadratic, growth in the number of operations with the growth of N.

## Conclusions

- Usage of the fast matrix multiplication method (the Dourbal method) significantly reduces the computational complexity of numerical solutions for partial differential equations. The number of addition and multiplication operations grows linearly with vector length N.
- The number of addition and multiplication operations is proportional to the required solution accuracy.
- The Dourbal method should be used where fast solutions of differential equations are needed.

## References

[1] Pavel Dourbal. " Synthesis of fast multiplication algorithms for arbitrary tensors." arXiv:1602.07008.

[2] Pavel Dourbal. " Method and apparatus for fast digital filtering and signal processing." US Patent Application # 14/748541, June 24, 2015.

[3] Saad Y., “Iterative Methods for Sparse Linear Systems”, N.Y.: PWS Publ., 1996. 547 p.

[4] Pan V. Ya., “Strassen's algorithm is not optimal — trilinear technique of aggregating uniting and canceling for constructing fast algorithms for matrix operations”, Proc. 19th Annual Symposium on Foundations of Computer Science, Ann Arbor, Mich., 1978.

[5] Bini D., Capovani M., Lotti G., Romani F., “$O(n^{2.7799})$ complexity for approximate matrix multiplication”. — Inform. Process. Lett., 1979.

[6] Schonhage A., “Partial and total matrix multiplication”, SIAM J. Comput., 1981.

[7] Don Coppersmith and Shmuel Winograd, “Matrix multiplication via arithmetic progressions”, Journal of Symbolic Computation, 9 p. 251–280, 1990.

[8] Virginia Vassilevska, Ryan Williams: Finding, Minimizing, and Counting Weighted Subgraphs. SIAM J. Comput. 42(3): 831-854 (2013).

[9] Louis A. Hageman, David M. Young, “Applied Iterative Methods” Academic Press 1981.